\documentclass[10pt,a4paper]{article}
\usepackage[latin1]{inputenc}
\usepackage[english]{babel}
\usepackage{graphics}
\usepackage{graphpap}
\usepackage{amssymb}
\usepackage{amsmath}
\newfont{\cc}{cmbsy10 scaled 1000}

\newcommand{\C}{{\mathbb C}}

\newcommand{\nbd}{neighbourhood }

\newcommand{\LL}{{\cal L}}
\newcommand{\Ln}{{\mathbb L}(n)}
\newcommand{\Lnn}[1]{{\mathbb L}^#1 (n)}
\newcommand{\LLn}{{\mathbb L}({\cal L})}
\newcommand{\LLt}{{\mathbb L}^1({\cal L})}

\newcommand{\M}{{\mathbb M}}

\newcommand{\N}{{\mathbb N}}
\newcommand{\R}{{\mathbb R}}
\newcommand{\Ss}{{\mathbb S}}
\newcommand{\T}{{\mathcal T}}
\newcommand{\Z}{{\mathbb Z}}

\newcommand{\ptir}{\discretionary{.}{}{.\kern.35em---\kern.7em}}
\newcommand{\cqfd}{\hfill\rule{5pt}{5pt}\vskip5pt}
\newcommand{\proof}{{\sl Proof}\ptir}
\newcommand{\sgn}{\rm sgn\,}
\newtheorem{Theo}{Theorem}[section]

\newtheorem{lem}{Lemma}[section]
\newtheorem{prop}{Proposition}[section]
\newtheorem{rem}{Remark}[section]
\newtheorem{cor}{Corollary}[section]
\newtheorem{definition}{Definition}[section]

\begin{document}
\selectlanguage{english}
\begin{titlepage}
\title
{A topological definition of the Maslov bundle}

\author
{{\Large Colette Anné}\\
anne@math.univ-nantes.fr\\
Laboratoire de Math\'ematiques Jean Leray\\
Universit\'{e} de Nantes, BP 92208 \\ 
44322 Nantes-Cedex 03, France
}
\date{\today}
\maketitle
\begin{abstract}We give a definition of the Maslov fibre bundle for a
lagrangian submanifold of the cotangent bundle of a smooth manifold. This
definition generalizes the definition given, in homotopic terms, by Arnol'd 
for lagrangian submanifolds of $T^\ast\R^n$. We show that our definition 
coincides with the one of Hörmander in his works about Fourier Integral
Operators.\\
MSClass : 81Q20, 35S30, 47G30.

\end{abstract}

\tableofcontents
\end{titlepage}

\section{Introduction}
The Maslov index appears as the phase term when one tries to define the symbol of a
Fourier Integral Operator (FIO). This symbol is then defined as a section of the 
Maslov bundle contructed on a lagrangian submanifold of $T^\ast X$. In his historical
paper [\ref{Ho2}], Hörmander proposes a construction of this bundle in terms of cocycles 
and tries to make the links with the strictly topological presentation (representation 
of the fundamental group) proposed by Arnol'd [\ref{arnold}], originally in an appendix of 
the book of Maslov [\ref{maslov}]. This link is established only for the lagrangian submanifolds 
of $T^\ast\R^n$. I propose in this work a new construction (\ref{def}) for the lagrangian 
submanifolds of $T^\ast X$, $X$ a smooth manifold, based on a definition of the Maslov index (\ref{def0}) 
which generalize the one of Arnol'd, and satisfies the cocycles conditions of Hörmander. These 
correspondances are established in the sections 2 and 3.

{\sl aknowledge}\ptir To François Laudenbach and his interest in this work.
\subsection{Arnol'd's definition of the Maslov index}\label{1}
Recall first the construction of Arnol'd [\ref{arnold}]. The space $T^\ast\R^n$
has a symplectic structure by the standard symplectic form
$$\omega=\sum_{j=1}^{j=n}d\xi_j\wedge dx_j.
$$
Let $\Ln$ be the Grassmannian manifold of the Lagrangian subspaces of $T^\ast\R^n$;
we identify $\Ln=U(n)/O(n).$ The map $Det^2$ is well defined on $\Ln.$ 
it is showed in [\ref{arnold}] that every path $\gamma :\ \Ss^1\to \Ln$ such that
$Det^2\circ\gamma\ :\,\Ss^1\to\Ss^1$ is a generator of $\Pi_1(\Ss^1),$ gives a
generator of $\Pi_1(\Ln).$ It follows that  $\Pi_1(\Ln)\simeq\Z$ and that the cocycle
$\mu_0$ defined by
$$\forall\gamma\in\Pi_1(\Ln)\quad\mu_0(\gamma)=\hbox{Degree }(Det^2\circ\gamma)$$
is a generator of the group $H^1(\Ln)\simeq\Z.$
It is then possible to define a {\em  Maslov bundle} ${\mathbb M}(n)$ on $\Ln$ by the
representation $\exp(i\frac{\pi}{2}\mu_0)=i^{\mu_0}$ of $\Pi_1(\Ln).$ It is
a flat bundle with torsion because ${\mathbb M}(n)^{\otimes 4}$ is trivial.

Now the Maslov bundle of a submanifold ${\LL}$ of $T^\ast\R^n$ is the pullback
of ${\mathbb M}(n)$ by the natural map
\begin{eqnarray*}
\varphi_n : \LL&\to& \Ln\\
\nu &\mapsto & T_\nu\LL.
\end{eqnarray*}
Arnol'd precisely shows that $\mu={\varphi_n}^\ast\mu_0$ is the Maslov index of $\LL.$
One can write
\begin{eqnarray}
\mu : \Pi_1(\LL)&\to& \Z\nonumber\\
\lbrack\gamma\rbrack &\mapsto & <\mu_0,\varphi_n\circ\gamma>=\hbox{Degree }(Det^2\circ\varphi_n\circ\gamma).
\end{eqnarray}
We have to take care of the structural group of this bundle. As a $U(1)-$bundle
it is always trivial. But it is concidered as a $\Z_4=\{1,i,-1,-i\}$-bundle.
In fact one can see, using the expression of the Maslov cocycle $\sigma_{jk}$ given by [\ref{Ho2}] 
(3.2.15) that the Chern classes of this bundle are null but $\sigma_{jk}$ can not be writen in general
as the coboundary of a {\em constant} cochain.

We recall now the theorem of symplectic reduction as it is presented in
[\ref{GS}] Proposition 3.2. p.132 .
\begin{prop}[Guillemin, Sternberg]\ptir\label{reduc} Let $\Delta$ be an isotropic subspace of dimension
m in $T^\ast\R^{(n+m)}.$ Define 
$S_\Delta=\{\lambda\in{\mathbb L}(n+m)/\;\lambda\supset\Delta\}.$ 
Then $S_\Delta$ is a submanifold of
${\mathbb L}(n+m)$ of codimension $(n+m)$, if we define $\rho$ to be the map 
\begin{eqnarray*}{\mathbb L}(n+m)&\stackrel{\rho}{\to}&{\mathbb L}(n)\\
\lambda&\mapsto&\lambda\cap \Delta^\omega /\lambda\cap \Delta
\end{eqnarray*}
($\Delta^\omega$ is the orthogonal of $\Delta$ for the canonical symplectic form   
$\omega$), then the map $\rho,$ which is continue on the all ${\mathbb L}(n+m),$ is
smooth in restriction to ${\mathbb L}(n+m)-S_\Delta$ and defines on this space a
fibre structure with base 
${\mathbb L}(n)$ and fibre $\R^{(n+m)}$. 

Moreover the image by $\rho$ of the generator of 
$\Pi_1({\mathbb L}(n+m))$  is a generator of $\Pi_1({\mathbb L}(n))$. 
\end{prop}

\subsection{Hörmander's definition of the Maslov bundle}\label{horm}
Let $X$ be a smooth manifold, then $T^\ast X\stackrel{\pi_0}{\to}X$ is endowed with a canonical
symplectic structure by $\omega=d \xi\wedge d x$. Let $\LL$ be a lagrangian (homogeneous)
submanifold of $T^\ast X$. Hörmander, in [\ref{Ho2}] p.155, defines the Maslov bundle of $\LL$ 
by its sections.

A Lagrangian manifold owns an atlas such that the cards $(C_\phi,D_\phi)$ are defined by 
non degenerated phase functions $\phi$ defined on  $U\times \R^N$ $U$ open in a domain diffeomorphic 
to a ball of a card of $X$ and
\begin{eqnarray*}
C_\phi=\Big\{(x,\theta);\,\phi'_\theta(x,\theta)=0\Big\}&\stackrel{D_\phi}{\longrightarrow}&
\LL_\phi\subset\LL\\
(x,\theta)&\longmapsto&(x,\phi'_x(x,\theta)).
\end{eqnarray*}
For the function $\phi$, to be non degenerate means that $\phi'_\theta$ is a submersion and thus
$C_\phi$ is a submanifold and $D_\phi$ an immersion.
\\
A section is then given by a family of functions
$$z_\phi : C_\phi\to \C
$$
satisfying the change of cards formulae~:
\begin{eqnarray}\label{compatible}
z_{\tilde{\phi}}=\exp i\frac{\pi}{4}\Big(\hbox{sgn}\phi''_{\theta\theta}-
\hbox{sgn}\tilde{\phi}''_{\tilde\theta\tilde\theta}\Big)z_\phi.
\end{eqnarray}
In fact $(\hbox{sgn}\phi''_{\theta\theta}-\hbox{sgn}\tilde{\phi}''_{\tilde\theta\tilde\theta})$ 
is even (see below, proposition \ref{sign}) and we have indeed constructed by this way a $\Z_4-$bundle.

\subsection{Definition of the Maslov  index and results}
In the same situation as before, we can construct on any lagrangian submanifold
$\LL$ of $T^\ast X$ (and in fact on all $T^\ast X$) the following fibre bundle
\[\begin{array}{ccc}
\Ln&\stackrel{i}{\longrightarrow}& \LLn\\
&& \pi\Big\downarrow\\
&&\LL \\
\end{array}\]
of the lagrangian subspaces of $T_\nu(T^\ast X), \nu\in \LL.$

This bundle has two natural sections~:
$$\lambda(\nu)=T_\nu(\LL),\;\hbox{and }\lambda_0(\nu)=\hbox{vert}(T_\nu(T^\ast X))
$$
defined by the tangent to $\LL$ and the tangent to the vertical $T_{\pi_0(\nu)}^\ast X.$

To a fibre bundle is associated a long exact sequence of homotopy groups, here~:
$$...\Pi_2(\LL)\to\Pi_1(\Ln)\stackrel{i_\ast}{\to}\Pi_1(\LLn)\stackrel{\pi_\ast}{\to}\Pi_1(\LL)
\to\Pi_0(\Ln)=0.
$$

But our fibre bundle possesses a section (two in fact), as a consequence the maps
$\Pi_k(\LLn)\stackrel{\pi_\ast}{\to}\Pi_k(\LL)$
are onto and the maps $\Pi_{k+1}(\LL)\to\Pi_k(\Ln)$ are null~; this gives a split exact sequence 
$$0\to\Pi_1(\Ln)\stackrel{i_\ast}{\to}\Pi_1(\LLn)\stackrel{\pi_\ast}{\to}\Pi_1(\LL)
\to 0.
$$
Take a base point $\nu_0\in\LL$ and fix a path $\sigma$ from
$\lambda(\nu_0)$ to $\lambda_0(\nu_0)$ lying in the fibre ${\LLn}_{\nu_0}$. 
For $\gamma\in\Pi_1(\LL)$ 
we denote ${{\lambda_0}^\sigma}_\ast(\gamma)$ the composition of  
$\sigma,\, {\lambda_0}_\ast\gamma$ and finaly $\sigma^{-1}$ (we use here the conventions of 
writing of [\ref{Husem}]).

Then
$\forall \gamma\in\Pi_1(\LL),\;\pi_\ast\Big({\lambda}_\ast\gamma\ast({{\lambda_0}^\sigma}_\ast(\gamma^{-1}))
\Big)=0$ and ${\lambda}_\ast\gamma\ast({{\lambda_0}^\sigma}_\ast(\gamma^{-1}))$ is in $\Pi_1(\Ln).$ 
Let us take the
\begin{definition}\label{def0}\ptir The Maslov index of $\LL$ is the map $\mu$ :
$$\forall \gamma\in\Pi_1(\LL),\;\mu(\gamma)=
\mu_0\Big({\lambda}_\ast\gamma\ast{{\lambda_0}^\sigma}_\ast(\gamma^{-1})\Big).
$$
\end{definition}

\begin{prop}\label{mu}\ptir This definition does not depend on the path $\sigma$ that we have
chosen to joint $\lambda(\nu_0)$ to $\lambda_0(\nu_0)$~; moreover $\mu$ is a morphism of group,
that is~: $\mu\in H^1(\LL,\Z)$.
\end{prop}

First remark~: in the case where $X=\R^n$ the fibre bundle $\LLn$ can be trivialized in such a way that
the section $\lambda_0$ is constant. In this case our definition coincide with the one of [\ref{arnold}]. 
A natural consequence of the proposition is the following definition~:

\begin{definition}\label{def}\ptir
The Maslov bundle $\M(\LL)$ over $\LL$ is defined as in section \ref{1} by the representation 
$\exp(i\frac{\pi}{2}\mu)=i^{\mu}$ of $\Pi_1(\LL)$ in $\C$.
\end{definition}
This means that the sections of the bundle are identified with functions $f$ on the
universal cover of $\LL$ with compex values and satisfying the relation~: 
\begin{equation}\label{equivariant}
\forall\gamma\in\Pi_1(\LL),\quad f(x.\gamma)=i^{-\mu(\gamma)}f(x),
\end{equation}
like in [\ref{amcharb}] formula (2.19).

%\begin{cor}\ptir Lorsque ceci a un sens, {\it i.e.} génériquement sur $\LL$, 
%$\mu$ est la duale de Poincaré de la sous-variété de $\LL$ des points où les 
%lagrangiens $\lambda$ et $\lambda_0$ ne sont pas transverses.
%\end{cor}

\begin{Theo}\label{defhorm}\ptir The sections of the Maslov  bundle of a Lagrangian (homogeneous) 
submanifold as defined by the definition \ref{def} satisfy the gluing conditions of Hörmander, it means
that our definition coincides with the one of Hörmander.
\end{Theo}

\section{Study of the index $\mu$.}

\subsection{The index $\mu_0$ on $\Ln$ is also an intersection number.} 
For $\alpha\in\Ln$ et $k\in\N$ one defines
$\Lnn{k}(\alpha)=\{\beta\in\Ln ; \; \hbox{dim } \alpha\cap\beta = k\}.$ 
Since [\ref{arnold}] we know that $\Lnn{k}(\alpha)$ is an open submanifold of codimension
$\frac{k(k+1)}{2},$ in particular $\overline{\Lnn{1}(\alpha)}$ is an oriented cycle of codimension 1 and
his intersection number coincides with $\mu_0$.

\subsection{Proof of the proposition \ref{mu}.}
It is a consequence of the two following lemmas. Provide $\LLn$ with a connection of $U(n)$-bundle.
Indeed any symplectic manifold $(M,\omega)$, like $T^*X$, can be provided with an almost complex 
structure $J$ which is compatible with the symplectic structure(see [\ref{matthias}] 
p.102), it means such that $g(X,Y)=\omega(JX,Y)$ is a riemannian metric. By this way the tangent
bundle of $M$ is provided with an hermitian form $g_\C=g+i\omega$, and its structural group
restricts to $U(n)$ it is also the case for the grassmannian of Lagrangians or its restriction to
a submanifold.

We will denote by $\tau(\gamma)_{x\to y}$ the parallel transport for this connection from $\LLn_{x}$ to $\LLn_{y}$ 
along the path $\gamma$  joining  $x$ to $y$ in $\LL.$  

Let's now $\gamma :\ \Ss^1\to \LL$ be a closed path such that $\gamma(0)=\nu_0$, we define
$\lambda(t)=\lambda_\ast(\gamma)(t)$ and in the same way  
${\lambda_0}^{-1}(t)={\lambda_0}_\ast(\gamma^{-1})(t).$

If, as before, $\sigma$ is a path from $\lambda(0)$ to
$\lambda_0(0)$ in the fibre $\LLn_{\gamma(0)}$~; then the path of $\LLn$ : 
$\lambda\ast\sigma\ast{\lambda_0}^{-1}\ast{\sigma}^{-1}$ is homotopic to
a path in the fibre, we have to calculate the Maslov index $\mu_0$ of this last one. For this we use 
the parallel transport along $\gamma$ to deform $\lambda\ast\sigma\ast{\lambda_0}^{-1}.$

\begin{definition}\ptir For $t\in[0,1]$ let's $\sigma_{t}$ denote the path included in the fibre 
$\LLn_{\gamma(t)}$ joining $\lambda(t)$ to $\lambda_0(t)$ and obtained by the parallel transport
of $\lambda_{|[t,1]}\ast\sigma\ast({{\lambda_0}}_{|[t,1]})^{-1}.$
\end{definition}
\begin{center}
\scalebox{.65}{\includegraphics{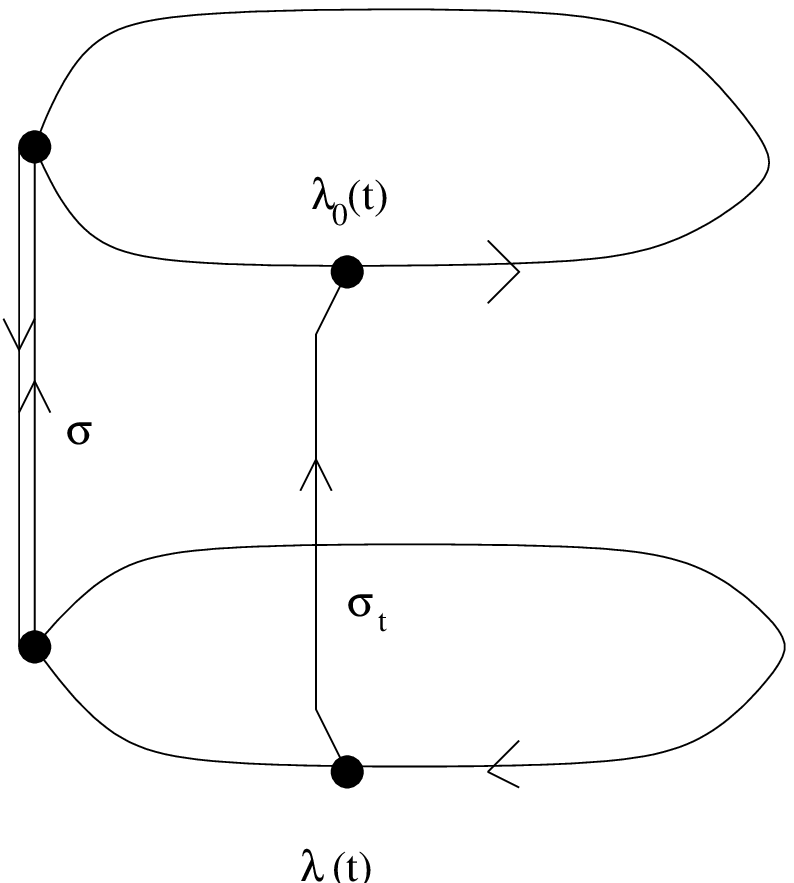}}
\end{center}

This path has three distinct parts : first 
$\tilde\lambda(t,s)=\tau(\gamma^{-1})_{\gamma(s)\to\gamma(t)}\lambda(s)$ then \break
$\tilde\sigma(t,s)=\tau(\gamma^{-1})_{\gamma(1)\to\gamma(t)}\sigma(s)$ and finally
${\tilde\lambda_0}^{-1}(t,s)=\tau(\gamma^{-1})_{\gamma(s)\to\gamma(t)}({\lambda_0}^{-1}(t)).$

By the definition (\ref{def})
$$\mu(\gamma)=\mu_0(\sigma_{0}\ast\sigma^{-1}).$$

\begin{lem}\ptir This definition does not depend on the path $\sigma$ chosen to link
$\lambda(0)$ to $\lambda_0(0)$ staying in the fibre above $\gamma(0).$
\end{lem}

The index $\mu_0$ is defined on the free homotopy group so
$$\mu_0(\sigma_{0}\ast\sigma^{-1})=\mu_0(\sigma^{-1}\ast\sigma_{0})=
\mu_0(\sigma^{-1}\ast\tilde\lambda\ast\tilde\sigma\ast{\tilde\lambda_0}^{-1})
$$
if, here, $\tilde\lambda(s)=\tilde\lambda(0,s)$ and the same notations for
$\lambda_0$ and $\sigma.$

If $\sigma'$ is an other path from $\lambda(0)$ to $\lambda_0(0)$, then
by the preceding remark and the fact that $\mu_0$ is a morphism of group, one has~: 
\begin{eqnarray*}
\mu_0({\sigma'_0}\ast{\sigma'}^{-1})-\mu_0(\sigma_{0}\ast\sigma^{-1})
&=&\mu_0({\sigma'}^{-1}\ast{\sigma'_0})-\mu_0(\sigma^{-1}\ast\sigma_{0})=\\
\mu_0({\sigma'}^{-1}\ast{\sigma'_0})+\mu_0({\sigma_{0}}^{-1}\ast\sigma)
&=&\mu_0({\sigma'}^{-1}\ast{\sigma'_0}\ast{\sigma_{0}}^{-1}\ast\sigma)=\\
\mu_0({\sigma'}^{-1}\ast\tilde\lambda\ast\tilde\sigma'\ast{\tilde\lambda_0}^{-1}
\ast({\tilde\lambda_0}^{-1})^{-1}\ast
{\tilde\sigma}^{-1}\ast{\tilde\lambda}^{-1}\ast\sigma)&=&
\mu_0(\sigma'^{-1}\ast\tilde\lambda\ast\tilde{\sigma'}\ast{\tilde\sigma}^{-1}
\ast{\tilde\lambda}^{-1}\ast\sigma)=\\
\mu_0(\sigma\ast{\sigma'}^{-1}\ast\tilde\lambda\ast\tilde{\sigma'}\ast{\tilde\sigma}^{-1}
\ast{\tilde\lambda}^{-1})&=&
\mu_0\Big((\sigma\ast{\sigma'}^{-1})\ast\tilde\lambda\ast(\tilde\sigma\ast{\tilde{\sigma'}}^{-1})^{-1}
\ast{\tilde\lambda}^{-1}\Big)=\\
\mu_0(\sigma\ast{\sigma'}^{-1})+\mu_0(\tilde\lambda\ast(\tilde\sigma\ast{\tilde{\sigma'}}^{-1})^{-1}
\ast{\tilde\lambda}^{-1})&=&\\
\mu_0(\sigma\ast{\sigma'}^{-1})+\mu_0({\tilde\lambda}^{-1}\ast\tilde\lambda\ast(\tilde\sigma\ast
{\tilde{\sigma'}}^{-1})^{-1})&=&\\
\mu_0(\sigma\ast{\sigma'}^{-1})+\mu_0((\tilde\sigma\ast
{\tilde{\sigma'}}^{-1})^{-1})
&=&\mu_0(\sigma\ast{\sigma'}^{-1})-\mu_0(\tilde\sigma\ast
{\tilde{\sigma'}}^{-1})=0
\end{eqnarray*}
because $\tilde\sigma\ast{\tilde{\sigma'}}^{-1}$ is the image of $\sigma\ast{\sigma'}^{-1}$ by
the parallel transport $\tau(\gamma)$ along $\gamma$ ; but $\tau(\gamma)\in U(n)$ preserves the Maslov 
index $\mu_0.$
\cqfd

\begin{lem}\ptir $\mu$ is a morphism of groups.
\end{lem}
Indeed, if $\alpha$ and $\beta$are two elements of $\Pi_1(\LL)$ it is suffisant to calculate
$\mu(\alpha)+\mu(\beta)$ beginning the first circle at $\tilde\sigma^{-1}(1)=
\tau(\alpha)\sigma(0)$ and applying $\tau(\alpha)$ to the second circle which was chosen to begin 
at $\sigma(0).$\cqfd

\section{Links with the definition of H\"ormander}
To make the link of this definition with signature terms of the formula in [\ref{Ho2}] we
follow the calculation from [\ref{duister}].
\subsection{Maslov's index in term of signature.} Let $\gamma\in\Lnn{k}(\alpha)$ and 
$\beta\in\Lnn{0}(\alpha)\cap\Lnn{0}(\gamma)$. Then 
$\alpha$ and $\beta$ are transversal and $\gamma$ can be presented as a graph~: there exists a unique 
linear map
$C:\alpha\to\beta$ such that $\gamma=\{(x,Cx), x\in\alpha\}.$
[\ref{duister}] p. 181,defines a quadratic form in $\alpha$ by~:
\begin{equation}
\label{dui}Q(\alpha,\beta;\gamma)=\omega(C.,.) \in{\cal Q}(\alpha).
\end{equation}
One sees easily that $\ker Q(\alpha,\beta;\gamma)=\ker C=\alpha\cap\gamma.$
and if we choose a basis on $\alpha$ such that $Q(\alpha,\beta;\gamma)$ has the form 
$\left|\begin{array}{cc}B_0&0\\0&0\end{array}\right|,$ the null part corresponds to
$\alpha\cap\gamma.$

Let now $\gamma(t)$ be a path in $\Lnn{0}(\beta)$ such that $\gamma(0)=\gamma.$ 
The goal of the following calculations is to control the jump of the signature
of the quadratic form $Q(\alpha,\beta;\gamma(t))$ in the neighbourhood of $t=0.$

\begin{prop}\label{duist1}\ptir Let $\gamma(t)$ be a path in $\Lnn{0}(\beta)$ such that
$\gamma(0)=\gamma$. If
$$Q(\alpha,\beta;\gamma(t))=\left|\begin{array}{cc}B(t)&C(t)\\C^t(t)&D(t)\end{array}\right|
$$
with $D(t)$ in  $\alpha\cap\gamma.$ Then, if $D'(t)$ is invertible in the neighbourhood of 0,
there exists $\varepsilon>0$ such that
$$\forall t,\;0<t<\varepsilon\;\;
\sgn Q(\alpha,\beta;\gamma(t))-\sgn Q(\alpha,\beta;\gamma(-t))= 
2\sgn D'(0).
$$
\end{prop}

\proof
We know that $B(t)$ is invertible and $C(t),\, D(t)$ are small.  The identity
\begin{equation}\label{diag}
\left|\begin{array}{cc}B&C\\C^t&D\end{array}\right|=
\left|\begin{array}{cc}1&0\\C^tB^{-1}&1\end{array}\right|.
\left|\begin{array}{cc}B&0\\0&(D-C^tB^{-1}C)\end{array}\right|.
\left|\begin{array}{cc}1&B^{-1}C\\0&1\end{array}\right|
\end{equation}
gives sgn $Q(\alpha,\beta;\gamma(t))=$sgn$(B(t))+$sgn$(D(t)-C(t)^tB(t)^{-1}C(t)).$
When $t$ is small sgn $B(t)=$ sgn $Q(\alpha,\beta;\gamma)$ and 
$\hbox{sgn}\Big(D(t)-C(t)^tB(t)^{-1}C(t)\Big)=\hbox{sgn}(t)\;\hbox{sgn}(D'(0))
$
by the mean value theorem.
\cqfd

Now if $\gamma$ is a path which cross  {\em transversally} $\Lnn{1}(\alpha)$
at $\gamma(0)$ then the assumption on $D'$ is satisfied.

\begin{Theo}\label{duist2}\ptir Let $\alpha\in\Ln$ and $\gamma$ a closed path in $\Ln$ which cross
$\Lnn{1}(\alpha)$ transversally, then for all $\beta\in\Ln$ transversal to $\alpha$ and to
$\gamma(t)$ one has 
$$\mu_0(\gamma)=\frac{1}{2}\sum_{t,\gamma(t)\in\Lnn{1}(\alpha)}
\Big(\sgn Q(\alpha,\beta;\gamma(t^+))-\sgn Q(\alpha,\beta;\gamma(t^-))  \Big).
$$
\end{Theo}

Indeed, in this case $T_\gamma\Ln/T_\gamma\Lnn{1}(\alpha)\sim S^2(\alpha\cap\gamma)$ which is oriented
by the positive-definite quadratic forms and $\sgn D'(0)=\pm 1,$ we use then the previous formula.

\begin{rem}\ptir This formula allows to define index of path not necessarely closed, 
see [\ref{salamon}].
\end{rem}

\subsection
{Hörmander's index.} Let $\alpha, \,\beta,\,\beta'$ be three elements of
$\Ln$ such that $\beta,\,\beta'\in \Lnn{0}(\alpha).$ For any path $\sigma$ joining
$\beta$ to $\beta'$ one defines 
$$[\sigma,\alpha]=\mu_0(\hat\sigma)
$$
where $\hat\sigma$ is the closed path obtained from $\sigma$ by linking its endpoints staying in
$\Lnn{0}(\alpha):$
$$\hat\sigma=\sigma*\sigma_\alpha\;\hbox{and}\;\sigma_\alpha\subset\Lnn{0}(\alpha). 
$$
The theorem (\ref{duist2}) shows that $[\sigma,\alpha]$ does not depend on the way
$\sigma$ is closed staying in $\Lnn{0}(\alpha).$
Let now $\alpha'$  be a point in $\Lnn{0}(\beta)\cap\Lnn{0}(\beta').$ The {\em index of Hörmander}
is the number
$$s(\alpha,\alpha';\beta,\beta')=[\sigma,\alpha']-[\sigma,\alpha]=
\mu_0(\sigma*\sigma_{\alpha'}*(\sigma*\sigma_\alpha)^{-1})=\mu_0(\sigma_{\alpha'}*{\sigma_\alpha}^{-1})
$$
because the calculation of $\mu_0$ does not depend on the base point in $\Ss^1.$

This index depends only on the four points in $\Ln$ and not on the paths~:

\begin{prop}\label{horm1}\ptir Let $\beta,\;\beta'\in\Lnn{0}(\alpha)\cap\Lnn{0}(\alpha')$ then
$$s(\alpha,\alpha';\beta,\beta')=\frac{1}{2}\Big(\sgn Q(\alpha,\beta';\alpha')-
\sgn Q(\alpha,\beta;\alpha')\Big).
$$
\end{prop}
Indeed, first suppose that $\alpha$ and $\alpha'$ are transversal ; the theorem (\ref{duist2})
can be applied and also the proposition (\ref{duist1}) ; this gives
$$s(\alpha,\alpha';\beta,\beta')=\frac{1}{2}\Big(\sgn Q(\alpha,\alpha';\beta)-\sgn  
Q(\alpha,\alpha';\beta')\Big).
$$
On the other hand $\beta\in\Lnn{0}(\alpha)$ can be writen as the graph of
$C\in \hbox{End}(\alpha,\alpha')$ and so $Q(\alpha,\alpha';\beta)=\omega(C.,.).$ But 
also $\alpha'$ is the graph of $D\in\hbox{End}(\alpha,\beta)$ with
$\forall x\in \alpha,\quad D(x)=-\big(x+C(x)\big),$
then $Q(\alpha,\beta;\alpha')=\omega(D.,.)=-\omega(C.,.)=-Q(\alpha,\alpha';\beta).$
As a consequence
$$s(\alpha,\alpha';\beta,\beta')=\frac{1}{2}\Big(\hbox{ sgn }Q(\alpha,\beta';\alpha')-\hbox{ sgn } 
Q(\alpha,\beta;\alpha')\Big).
$$

This formula can be generalized by the symplectic reduction (\ref{reduc}).
\cqfd

Let us recall finally the

\begin{prop}\label{horm2}\ptir Let $\alpha,\alpha',\beta,\beta'$ be four points in $\Ln$
such that $\beta$ and $\beta'$ are in $\Lnn{0}(\alpha)\cap\Lnn{0}(\alpha')$ then
$$s(\alpha,\alpha';\beta,\beta')=-s(\alpha',\alpha;\beta,\beta')=-s(\alpha,\alpha';\beta',\beta)
=-s(\beta,\beta';\alpha,\alpha').
$$
\end{prop}

Only the third equality is not obvious. It can be shown by the formula of proposition
\ref{horm1}. Choose symplectic coordinates $(x,\xi)$ such that $\alpha=\{x=0\}$ and 
$\beta=\{\xi=0\}.$
By the transversality hypothesis there exist homomorphisms $A$ and $B$ such that
$$\alpha'=\{x=A\xi\}\quad\beta'=\{\xi=Bx\}.
$$
If $\alpha'$ is the graph of $A'\in$ Hom$(\alpha,\beta'),$ then for all
$\xi\in\alpha$ we must find $\xi'\in\alpha$ and $x\in\beta$ with
$$A'\xi=(x,Bx)\;\hbox{and}\;(A\xi',\xi')=(x,Bx+\xi).
$$
This gives $x=A\xi'$ and $\xi'=Bx+\xi=BA\xi'+\xi$ so $\xi'=(1-BA)^{-1}\xi$ and
$$A'\xi=\big(A(1-BA)^{-1}\xi,(1-BA)^{-1}\xi-\xi\big).$$

We remark that $(1-BA)$ is indeed invertible : if $\xi\in\ker(1-BA)$ then
$(A\xi,\xi)=(A\xi,BA\xi)\in\alpha'\cap\beta'=\{0\}$ so $\xi=0.$

Therefore by the proposition (\ref{horm1})
$$2s(\alpha,\alpha';\beta,\beta')=\hbox{ sgn }\omega(A(1-BA)^{-1}.,.)-\hbox{ sgn }\omega(A.,.)=
\hbox{ sgn }\left|\begin{array}{cc}A&0\\0&-A(1-BA)^{-1}\end{array}\right|.
$$
Suppose now that $A$ is inversible then, because a symmetric matrix and its inverse have 
same signature~:
\begin{multline*}
\hbox{ sgn }\left|\begin{array}{cc}A&0\\0&-A(1-BA)^{-1}\end{array}\right|=
\hbox{ sgn }\left|\begin{array}{cc}A&0\\0&-(1-BA)A^{-1}\end{array}\right|=\\
=\hbox{ sgn }\left|\begin{array}{cc}A&0\\0&B-A^{-1}\end{array}\right|
=\hbox{ sgn }\left|\begin{array}{cc}A&1\\1&B\end{array}\right|
\end{multline*}
by formula (\ref{diag}). By the same calculus, and because $\omega$ is skewsymmetric,
one has~:
$$2s(\beta,\beta';\alpha,\alpha')=\hbox{ sgn }Q(\beta,\alpha';\beta')-\hbox{ sgn } 
Q(\beta,\alpha;\beta'))=-\hbox{ sgn }\left|\begin{array}{cc}B&1\\1&A\end{array}\right|.
$$
\cqfd

\subsection{Proof of theorem \ref{defhorm}} 
Following [\ref{Ho2}], we denote by $\T(\LL)\subset\LLn$ the set of the $\alpha\in\LLn$
transversal to $\lambda(\pi(\alpha))$ and to $\lambda_0(\pi(\alpha)).$
If $p:\T(\LL)\to\LL$ is the associated projection, then for all $\nu\in\LL$
$$p^{-1}(\nu)=\Lnn{0}(\lambda(\nu))\cap\Lnn{0}(\lambda_0(\nu)).$$
{\it n.b.} On the \nbd of points where the two Lagrangian are not transversal this map
is not a fibration.

\begin{lem}\label{continu}\ptir Let $\alpha:\ \Ss^1\to \T(\LL)$ satisfying $p\circ\alpha=\gamma$
and $\sigma$ be a path as before. The index 
$[\sigma_t,\alpha(t)]$ is constant in $t$. 
\end{lem}

Indeed the index is a continuous map~: let $t_0\in [0,1]$ and $\beta$ a path in the fibre over 
the point $\gamma(t_0)$ and linking $\lambda_0(t_0)$ to $\lambda(t_0)$ staying tranversal to
$\alpha(t_0)$; by definition $[\sigma_{t_0},\alpha(t_0)]=\mu_0(\sigma_{t_0}\ast\beta)$
but the property of transversality is open~: if we denote $\beta_t$ the path in the fibre over 
the point $\gamma(t)$ resulting of the parallel transport of 
${\lambda_0}_{|[t, t_0]}*\beta*{\lambda^{-1}}_{|[t, t_0]},$ then there exists $\varepsilon>0$ 
such that for all $|t-t_0|<\varepsilon$ one has $\beta_t$ is transversal to $\alpha(t).$
This parallel transport realizes an homotopy, so for all
$|t-t_0|<\varepsilon$ one has
$\mu_0(\sigma_{t_0}\ast\beta)=\mu_0(\sigma_{t}\ast\tilde\beta_t).$ 
\cqfd

\begin{cor}\ptir The induced fibrex bundle $p^*\M(\LL)$ is trivial.
\end{cor}

\proof We have to show that for all path $\alpha : \Ss^1\to \T(\LL)$  continuous, if we define
$\gamma=p\circ\alpha,$ then $\mu(\gamma)=0.$
To this goal take $\sigma$ as before, a path in the fibre over $\gamma(0)$ linking $\lambda(0)$ to 
$\lambda_0(0).$ Choose $\sigma$ transversal to $\alpha(1)$ and do the same constrution as before, then
$$[\sigma,\alpha(1)]=[\sigma_0,\alpha(0)]=0
$$
by the definition of $[\sigma,\alpha(1)]$ and lemma \ref{continu}.
But $\alpha(0)=\alpha(1)$ so 
$$\mu(\gamma)=\mu_0(\sigma_0*\sigma^{-1})=[\sigma_0,\alpha(1)]=0.
$$
\cqfd
\begin{cor}\ptir Let $s$ be a section of the Maslov bundle over $\LL,$ and $\gamma :\ \Ss^1\to \LL$ 
a closed path such that $\gamma(0)=\nu_0=\pi(\lambda_0).$ Let $\alpha:\ [0,1]\to \T(\LL)$ be a continuous path
satisfying $\gamma=p\circ\alpha.$ Then
$$p^\ast s(\alpha(1))=i^{s(\lambda_0(0),\lambda(0);\alpha(1),\alpha(0))}p^\ast s(\alpha(0)).
$$
\end{cor}
\proof  Let $\sigma$ be a path linking $\lambda(0)$ to $\lambda_0$  staying transversal to
$\alpha(1)$. By lemma (\ref{continu}), $[\sigma_0,\alpha(0)]$ $=[\sigma,\alpha(1)]=0$ and 
$$\mu(\gamma)=\mu_0(\sigma_0\ast\sigma^{-1})=[\sigma_0,\alpha(1)]=[\sigma_0,\alpha(1)]-[\sigma_0,\alpha(0)]
=s(\alpha(0),\alpha(1);\lambda(0),\lambda_0(0))
$$
and $s(\alpha(0),\alpha(1);\lambda,\lambda_0)=-s(\lambda_0,\lambda;\alpha(1),\alpha(0))$
by the proposition \ref{horm2}. Therefore 
$$-\mu(\gamma)=s(\lambda_0(0),\lambda(0);\alpha(1),\alpha(0)).
$$
This gives the result by the equivalent relation (\ref{equivariant}).
\cqfd

From these two corollaries one obtains 
\begin{cor}\ptir The sections of $\M(\LL)$ are identified with functions $f$ on $\T(\LL)$ 
satisfying the relation~: $\forall \alpha,\tilde{\alpha}\in\T(\LL)$
$$p(\alpha)=p(\tilde{\alpha})\Rightarrow f(\tilde{\alpha})=i^{s(\lambda_0,\lambda;\tilde{\alpha},\alpha)}
f(\alpha).
$$
\end{cor}
This result gives the gluing condition of Hörmander, in view of the theorem 3.3.3, [\ref{Ho2}] and finish
the proof of the theorem. For completness we recall this last step.
\begin{prop}\label{sign}\ptir The functions $f$ on $\T(\LL)$ which satisfy~: $\forall \alpha,
\tilde{\alpha}\in\T(\LL)$
$$p(\alpha)=p(\tilde{\alpha})\Rightarrow f(\tilde{\alpha})=i^{s(\lambda_0,\lambda;\tilde{\alpha},\alpha)}
f(\alpha).
$$
are the sections defined by the gluing conditions of the section \ref{horm}.
\end{prop}
\proof Let $\phi$ be a non degenerated phase function as in section \ref{horm} and
$\nu_0=(x_0,\xi_0)=(x_0,\phi'_{x}(x_0,\theta_0))$ a point in $\LL_\phi.$ For each 
$\alpha\in\T(\LL)$ such that $p(\alpha)=\nu_0,$ there exists a function $\psi$ defined on an open set $U$ 
such that the graph $L_\psi=\{(x,d\psi(x)),\,x\in U\}$ of the differential $d\psi$ intersect transversally 
$\LL_\phi$ at $\nu_0,$ one has $\xi_0=d\psi(x_0)$ and $T_{\nu_0} L_\psi=\alpha.$

Or equivalently one can say~: the following quadratic form defined on $\R^{n+N}$ by the
matrix
\begin{eqnarray}\label{Q}
Q_\psi=\left|\begin{array}{cc}\phi''_{xx}-\psi''_{xx}&\phi''_{x\theta}\\
                               \phi''_{\theta x} &\phi''_{\theta\theta}  
\end{array}\right| 
\end{eqnarray}
is non degenerated.

The restriction of this quadratic form to the tangent $W$ of $\LL_\phi$ at $\nu_0$  only depends on $\LL$ and 
$\psi$ (and not on $\phi$). Indeed $\phi$ defines a card in which
$$\lambda(\nu_0)=T_{\nu_0}(\LL)=\{(X,\phi''_{xx}X+\phi''_{x\theta}A);(X,A)\in\R^{n+N},\,
\phi''_{\theta x} X+\phi''_{\theta\theta} A=0\};
$$
if now $(X,A),(X',A')$ define two tangent vectors $V$ and $V'\in\lambda(\nu_0)$
\begin{eqnarray*}
Q_\psi\Big((X,A),(X',A')\Big)&=&<X,(\phi''_{xx}-\psi''_{xx})X'+\phi''_{x\theta}A'>\\
<-\psi''_{xx}X,X'>-<-X,\phi''_{xx}X'+\phi''_{x\theta}A'>
                             &=&Q\Big(\lambda(\nu_0),\alpha;\lambda_0(\nu_0)\Big)(V,V')\\
\end{eqnarray*}
by definition (\ref{dui}).
 More precisely $\alpha$ is transverse to the two lagrangians $\lambda(\nu_0)$ and $\lambda_0(\nu_0)$
so the vertical $\lambda_0(\nu_0)$ is the graph of an homomorphism $A_\psi$
from $\lambda(\nu_0)$ to $\alpha=T_{\nu_0} L_\psi$~: 
$$\forall (0,\Xi)\in\lambda_0(\nu_0),\;\exists(X,A)\hbox{unique such that } \Xi=\phi''_{xx}X+\phi''_{x\theta}A\,
\hbox{ et } \phi''_{\theta x} X+\phi''_{\theta\theta} A=0$$
because $Q_\psi$ is non degenerated, and one can write 
$$(0,\Xi)=(X,\phi''_{xx}X+\phi''_{x\theta}A)-(X,\psi''_{xx}X),
$$ 
it means that $A_\psi(X,\phi''_{xx}X+\phi''_{x\theta}A)=(-X,-\psi''_{xx}X).$

We see now that the orthogonal $W^{Q_\psi}$ of $W$ with respect to $Q_\psi$ 
is $\R^N=\{(0,A)\}$ and that 
${Q_\psi}_{|W^{Q_\psi}}=\phi''_{\theta\theta}.$
But the lemma \ref{sum} below gives 
$\sgn Q_\psi=\sgn {Q_\psi}_{|W}+\sgn {Q_\psi}_{|W^{Q_\psi}},$ so~:
\begin{eqnarray}\label{signature}
\sgn Q_\psi=\sgn Q(\lambda(\nu_0),\alpha;\lambda_0(\nu_0))+\sgn 
\phi''_{\theta\theta}.
\end{eqnarray}

Let now $z_\phi$ be a section in the sens of Hörmander. For any $\alpha\in\T(\LL), p(\alpha)=\nu_0,$
if $\phi$ and $\tilde\phi$ are two phase functions defining $\LL$ in a \nbd of $\nu_0$ and if $\psi$ 
is a  function on $X$ satisfying $\alpha=T_{\nu_0} L_\psi$, we denote by
$Q_\psi$ and $\tilde{Q}_\psi$ the respective quadratic forms defined by (\ref{Q}). Put
$$f(\alpha)=\exp(i\frac{\pi}{4}\sgn Q_\psi)z_\phi(\nu_0).$$
By the relation (\ref{signature}) one has 
$\sgn \phi''_{\theta\theta}-\sgn \tilde\phi''_{\tilde\theta\tilde\theta}=
\sgn Q_\psi-\sgn \tilde Q_\psi;
$
the compatibility condition \ref{compatible} gives then
$$\exp(i\frac{\pi}{4}\sgn Q_\psi)z_\phi(\nu_0)=
\exp(i\frac{\pi}{4}\sgn \tilde{Q}_\psi)z_{\tilde\phi}(\nu_0)
$$
and the function $f$ is well defined on $\T(\LL).$
On the other hand if $\tilde\alpha$ is an other point in $\T(\LL)$ such that
$p(\tilde\alpha)={\nu_0}$ and if $\tilde\psi$ is an adapted function, then
\begin{eqnarray*}
f(\tilde\alpha)&=&\exp(i\frac{\pi}{4}(\sgn \tilde Q_\psi-\sgn Q_\psi))f(\alpha)\\
&=&\exp\Big(i\frac{\pi}{4}\Big(\sgn Q(\lambda(\nu_0),\tilde\alpha;\lambda_0(\nu_0))-
\sgn Q(\lambda(\nu_0),\alpha;\lambda_0(\nu_0)\Big)\Big)f(\alpha)\\
&=&\exp\Big(i\frac{\pi}{2}s(\lambda(\nu_0),\lambda_0(\nu_0);\alpha,\tilde\alpha)\Big)f(\alpha)\\
&=&\exp\Big(i\frac{\pi}{2}s(\lambda_0(\nu_0),\lambda(\nu_0);\tilde\alpha,\alpha)\Big)f(\alpha)
\end{eqnarray*}
So it is a section of the Maslov bundle and the theorem \ref{defhorm} is proved.
\cqfd
\begin{lem}\ptir\label{sum} Let $Q$ be a  non degenerated quadratic form defined on
$\R^n,$ $V$ be a subspace of $\R^n$ and $V^Q$ its orthogonal for $Q$, then
$$\sgn Q=\sgn Q_{|V}+\sgn Q_{|V^Q}.
$$
\end{lem}
\proof This lemma can be showed using an  induction on dim$V\cap V^Q$. If  dim$V\cap V^Q=0$ there is 
nothing to do, if not let $v_1,\dots,v_k$ be a base of $V\cap V^Q.$ We complete this base with  
$v_{k+1},\dots,v_p$ to obtain a base of $V+ V^Q$. Because $Q$ is non degenerated there exists $w_1\in\R^n$
such that $Q(v_1,w_1)=1,$ and eventually after a modification with a linear combination of the $v_j$
one can suppose $Q(w_1)=0$ and $Q(v_j,w_1)=0$ for $j>1.$ One remarks that the signature of $Q$ in 
restriction to $\R v_1\oplus\R w_1$ is zero and applies the induction hypotheses to $(\R v_1\oplus\R w_1)^Q.$
\cqfd

\section{Topological comments}
Let's have a look to the exact sequence~:  
$0\to\Pi_1(\Ln)\stackrel{i_\ast}{\to}\Pi_1(\LLn)\stackrel{\pi_\ast}{\to}\Pi_1(\LL)
\to 0.$

The group $\Pi_1(\LLn)$ is the semidirect product of $\Pi_1(\Ln)$ and $\Pi_1(\LL).$ It means that
$\Pi_1(\LL)$ acts on $\Pi_1(\Ln)$ by conjugation. More precisely for all
$\gamma\in \Pi_1(\LL)$ let's define
\begin{eqnarray*}\rho_\gamma:\Pi_1(\Ln)&\to &\Pi_1(\Ln)\\
\sigma&\mapsto& \lambda_0(\gamma)\ast i_\ast(\sigma)\ast(\lambda_0(\gamma))^{-1}
\end{eqnarray*}
\begin{lem}This representation is trivial and $\Pi_1(\LLn)$ is in fact the direct product of
$\Pi_1(\Ln)$ and $\Pi_1(\LL).$
\end{lem}
\proof As was seen in paragraph 2, the parallel transport along  $\gamma$
defines an homotopy of $\lambda_0(\gamma)\ast i_\ast(\sigma)\ast(\lambda_0(\gamma))^{-1}$
to a path which can be writen $\tilde\lambda_0\ast\tilde\sigma\ast(\tilde\lambda_0)^{-1}$
where $\tilde\sigma$ is the image of $\sigma$ by $\tau(\gamma).$ But
$$\mu_0(\tilde\lambda_0\ast\tilde\sigma\ast(\tilde\lambda_0)^{-1})=
\mu_0((\tilde\lambda_0)^{-1}\ast\tilde\lambda_0\ast\tilde\sigma)=\mu_0(\tilde\sigma)=\mu_0(\sigma).
$$
As a consequence of the works of  Arnol'd recalled above, a generator of $\Pi_1(\Ln)$ is caracterized
by $\mu_0(\sigma)=1$.
\cqfd
\begin{Theo}\label{dual}\ptir Let $\LLt$ be the set of the points
$l\in\LL$ which are not transversal to $\lambda_0(\pi(l))$. It is an oriented cycle of $\LL$ 
of codimension 1~; if $m$ is its Poincaré dual form, then
$$\mu=\lambda^*m.
$$
\end{Theo}

\proof 
We keep the notations of paragraph 2. By choosing the starting point one can suppose that
{\em the two lagrangians $\lambda_0=\lambda_0(0)$ and $\lambda(0)$ are transversal.}
We will use a deformation of the path
$\tilde\lambda\ast\tilde\sigma\ast{\tilde\lambda_0}^{-1}$ joining $\lambda(0)$
to $\lambda_0(0).$ Recall that $\tilde\sigma(t)=\tau(\gamma)(\sigma(t)).$

There exists a (continuous) path $u(t)\in U(n)$ such that $u(0)=$I and
\[\forall t\in [0,1]\; \tilde\lambda_0(t)=u(t)(\lambda_0).
\]
But $\tilde\lambda_0(1)=\tau(\gamma)(\lambda_0),$ so $\tau(\gamma)$ and $u(1)$
differ by an element of $O(n)$~:
\[\exists a\in O(n)\; ;\; \tau(\gamma)=u(1)\circ a.
\]
Let's construct the following homotopy of $\tilde\lambda\ast\tilde\sigma\ast{\tilde\lambda_0}^{-1}$ 
by the concatenation of $u(st)^{-1}\tilde\lambda(t),$ next $u(s)^{-1}\tilde\sigma$ and finally the
inverse of $u(st)^{-1}\tilde\lambda_0(t).$ The end of this homotopy is a path, result of the 
concatenation of $\bar\lambda(t)=u(t)^{-1}\tilde\lambda(t)$ and $u(1)^{-1}\tilde\sigma=a\sigma$ 
because $u(t)^{-1}\tilde\lambda_0(t)=\lambda_0$ is a constant path.

We have now to calculate $\mu_0(\sigma^{-1}*\bar\lambda*a\sigma).$ Because $a\in O(n)$
\[Det^2(\sigma(t))=Det^2(a\sigma(t));
\]
$Det^2\circ\bar\lambda$ is a closed path even if $\bar\lambda$ is not, so
$\mu(\gamma)=\hbox{Degree }(Det^2\circ\bar\lambda).$

Considering the results of section 2.1, we have obtained
\begin{prop}
$\mu(\gamma)$ is the intersecting number of the submanifold $\overline{\Lnn{1}(\lambda_0)}$ and the
cycle obtained from $\bar\lambda$, by closing it with a path staying transversal to $\lambda_0$.
\end{prop}
Remark that $\bar\lambda(0)=\lambda(0)$ and  $\bar\lambda(1)=a\lambda(0)$ are both transversal to 
$\lambda_0.$
Let's now
$$\LLt=\Big\lbrace l\in \LLn\; ;\; \lambda_0(\pi(l))\cap l\neq \{0\}\Big\rbrace.
$$
It is a fibration above $\LL$ with fibre $\overline{\Lnn{1}(\lambda_0)}$, so it is an oriented cycle of
codimension 1 in $\LL.$ If $\lambda\circ\gamma$ cuts $\LLt$ transversally
at $\lambda\circ\gamma(t)$ then $\bar\lambda$ cuts transversally
$\overline{\Lnn{1}(\lambda_0)}$ at $\bar\lambda(t)$ and conversely. Moreover the transformations 
which permit to pass from $\lambda\circ\gamma$ to $\bar\lambda$ realise a continuous deformation 
of $\LLt$ to $\overline{\Lnn{1}(\lambda_0)}$ above $\gamma.$
This argument finishes the proof of the theorem \ref{dual}. \cqfd

\section{References}.
%\begin{thebibliography}{99}
\begin{enumerate}
\item\label{matthias} Aebischer B, Borer M \& al\ptir {\em Symplectic Geometry}, Birh\"auser, Basel, 1992.
\item\label{amcharb} Anné C., Charbonnel A-M\ptir {\em Bohr Sommerfeld conditions for several commuting
Hamiltonians}, preprint (Juillet 2002), ArXiv Math-Ph/0210026.
\item\label{arnold} Arnol'd V.I\ptir {\em Characteristic Class entering in Quantization Conditions.}
Funct. Anal. and its Appl. {\bf 1} (1967), 1-14.
\item\label{duister}Duistermaat J\ptir {\em Morse Index in Variational Calculus.} Adv. in Maths 
{\bf 21} (1976), 173--195. 
\item\label{DG} Duistermaat J., Guillemin V\ptir  {\em The spectrum of positive elliptic operators 
and periodic bicharacteristics.} Invent. Math. {\bf 29} (1975), 39-79.
\item\label{GS} Guillemin V., Sternberg S\ptir  {\em Geometric Asymptotics} Math. Surveys and 
Monograph n$^o$ 14, AMS (1990).
\item\label{Ho2} H\"{o}rmander L\ptir  {\em Fourier Integral operators.} 
Acta Math. {\bf 127} (1971), 79--183.
\item\label{Ho1} H\"{o}rmander L\ptir  {\em The Weyl calculus of pseudodifferential operators.} 
Comm. Pure Appl. Math. {\bf 32} (1979), 359-443.
\item\label{Ho3} H\"{o}rmander L\ptir  {\em The  Analysis  of  Linear  Partial
Differential Operators III.}, Springer, Berlin - Heidelberg  -  New  York, 1985.
\item\label{Ho4} H\"{o}rmander L\ptir  {\em The  Analysis  of  Linear  Partial
Differential Operators IV.}, Springer, Berlin - Heidelberg  -  New  York, 1985.
\item\label{Husem} Husemoller D\ptir {\em Fibre Bundles}, Springer, Berlin - Heidelberg  -  New  York, 1975.
\item\label{maslov} Maslov V.P\ptir {\em Théorie des perturbations et méthodes asymptotiques}, Dunod, Gauthier-Villars, Paris 1972.
\item\label{salamon} Robin J., Salamon D\ptir {\em The Maslov Index for Path.} Topology {\bf 32}
(1993), 827--844.
\end{enumerate}
%{thebibliography}
\end{document}